\documentclass[12pt]{amsart}
\usepackage{amsfonts}
\usepackage{amssymb}
\usepackage[mathscr]{eucal}
\usepackage{graphicx}

\newtheorem{theorem}[]{Theorem}
\newtheorem{lemma}[theorem]{Lemma}

\theoremstyle{definition}
\newtheorem*{definition}{Definition}
\newtheorem{example}[theorem]{Example}

\theoremstyle{remark}

\numberwithin{equation}{section}

\newcommand{\E}{\mathcal{E}}
\newcommand{\R}{\mathcal{R}}
\newcommand{\Pk}{\mathcal{P}_k}

\begin{document}

\title{Walks Along Braids and the Colored Jones Polynomial}

\author{Cody Armond}
\address{Mathematics Department\\
Louisiana State University\\
Baton Rouge, Louisiana}
\email{carmond@math.lsu.edu}

\date{}

\begin{abstract}
Using the Huynh and L\^{e} quantum determinant description of the colored Jones polynomial, we construct a new combinatorial description of the colored Jones polynomial in terms of walks along a braid. We then use this description to show that for a knot which is the closure of a positive braid, the first $N$ coefficients of the $N$-th colored Jones polynomial are trivial.
\end{abstract}

\maketitle

\section{Introduction}

The normalized colored Jones polynomial $J'_K(N)$ for a knot $K$ and positive integer $N \geq 2$ is a Laurent polynomial in the variable $q$, i.e. $J'_K(N) \in \R = \mathbb{Z}[q,q^{-1}]$. This is defined so that $J'_K(2)$ is the ordinary Jones polynomial, and $J'_U(N) = 1$ where $U$ is the unknot. In ~\cite{hl} Vu Huynh and Thang L\^{e} showed that the colored Jones polynomial of a knot $K$ can be described in terms of the inverse of the quantum determinant of an almost quantum matrix. This matrix comes from a braid $\beta$ whose closure is the knot $K$ and is formed using deformed Burau matrices. In ~\cite{j}, Vaughan Jones briefly mentioned a probabilistic interpretation of the  Burau representation as walks along the knot, and in ~\cite{ltw}, Xiao-Song Lin, Feng Tian, and Zhenghan Wang used this interpretation to generalize the Burau representation to tangles by using walks along the tangles. Then in ~\cite{lw} Lin and Wang used this to calculate the colored Jones polynomial and derive a new proof of the Melvin-Morton conjecture. In this paper, we will use a different, but similar, generalization of this description using Huynh and L\^{e}'s result to give a geometric interpretation of the colored Jones polynomial in terms of walks along the braid $\beta$.

This interpretation will then be used to show the following theorem, which basically states that, for knots which can be written as the closure of a positive braid, the first $N$ leading coefficients of $J'_K(N)$ are trivial:

\begin{theorem}
\label{positive}
 If $\beta$ is a positive braid whose closure is the knot $K$, then
$$J'_K(N)=\sum_{i=1}^N a_{N,i} q^{L_N+i-1} + \text{Higher terms}$$
where $a_{N,1} = 1$ and $a_{N,i} = 0$ for $2\leq i \leq N$.

Furthermore, $L_N=(N-1)(n-m+1)/2$ where $n$ is the number of crossings in $K$ and $m$ is the number of strands in the braid $\beta$.
\end{theorem}

In ~\cite{ck}, Abhijit Champanerkar and Ilya Kofman prove a result similar to Theorem ~\ref{positive} for a smaller class of knots, closures of positive braids with full twists. In their theorem, they show that this pattern of coefficients repeats multiple times using a different method.

\subsection{Plan of paper}
In section ~\ref{walks} walks are introduced and it is shown how they are used to calculate the colored Jones polynomial. In section ~\ref{qd} we discuss Huynh and L\^{e}'s description of the colored Jones polynomial as the inverse of a quantum determinant. This is then related to walks and is used to prove that the method described in section ~\ref{walks} gives us the colored Jones polynomial. Section ~\ref{posbraids} gives the proof of Theorem ~\ref{positive}.

\subsection{Acknowledgements}
I would like to thank Oliver Dasbach, whose guidance made this possible, and also Stavros Garafoulids and Thang L\^{e} for their hospitality during a visit to Georgia Tech University, and for helpful feedback.

\section{Walks}
\label{walks}
The braid group $B_m$ is the group of braids on $m$ strands. The product of two braids $\beta_1 \beta_2$ is the braid formed by putting the first braid $\beta_1$ on top of the second $\beta_2$. Let $\{\sigma_i\}_{1\leq i\leq m-1}$ be the standard generators of $B_m$.

The closure of a braid $\beta$ is the knot or link $\hat{\beta}$ formed by attaching the top of each strand to the bottom of each strand. Every knot can be realized as the closure of some braid. This braid is not unique, and neither is the word in the generators of the braid group which represents that braid. In this discussion we would like to fix such a word. The sequnce $\gamma$ which is described in the next paragraph gives us the means to do this.

For a sequence $\gamma = (\gamma_1,\gamma_2, \ldots, \gamma_k)$ of pairs $\gamma_j = (i_j,\epsilon_j)$, $1 \leq i_j \leq m-1$ and $\epsilon_j = \pm$ (we will use $\pm = \pm 1$ interchangably as it should be clear from the context), let $\beta = \beta(\gamma)$ be the braid
$$\beta := \sigma_{i_1}^{\epsilon_1}\sigma_{i_2}^{\epsilon_2}\ldots\sigma_{i_k}^{\epsilon_k}.$$
Fix a $\gamma$ so that the closure of $\beta(\gamma)$ is the knot $K$. Denote $\omega(\gamma) = \sum_{j=1}^k \epsilon_j$. This is the writhe of $\beta(\gamma)$.

\begin{definition}
 
We will define a path along the braid $\beta(\gamma)$ from $i$ to $j$ as follows:

Begining at the bottom of the $i$-th strand, follow the braid along a strand until you begin to cross over another strand. At this over crossing there is a choice to either continue along the strand or jump down to the strand below and continue following along the braid.  Continue to the top of the braid ending at the $j$-th strand.

\end{definition}

Each path is given a weight defined as follows:

At the $j$-th crossing, (i.e. the crossing corresponding to $\gamma_j$):
\begin{itemize}
\item If the path jumps down, assign $a_{j,\epsilon_j}$.
\item If the path follows the lower strand, assign $b_{j,\epsilon_j}$.
\item If the path follows the upper strand, assign $c_{j,\epsilon_j}$.
\end{itemize}

The weight of the path is the product of the weights of the crossings. Shortly we will give more meaning to these symbols, and then it will be immediate that for different values of $j$, they all commute with each other, so the order of the product can be taken arbitrarily. However, for the same value of $j$, these symbols are non-commutative, and we will see the relations between them.

\begin{figure}[htbp] %
   \centering
   \includegraphics[width=1.5in]{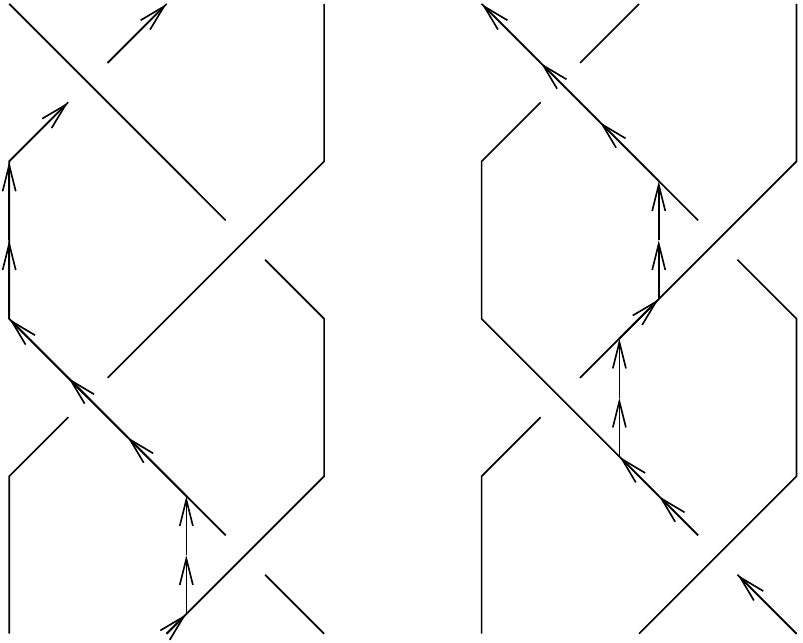} 
   \caption{A path from 2 to 2 and from 3 to 1}
   \label{paths}
\end{figure}

A walk $W$ along $\beta$ consists of a set $J \subset \{1,\ldots, m\}$, permutation $\pi$ of $J$, and a collection of paths where there is exactly one path from $j$ to $\pi(j)$ in the collection, for each $j \in J$.

The weight assigned to a walk is $(-1)(-q)^{|J| + \text{inv}(\pi)}$ times the product of the weights of the paths in the collection, where inv($\pi$) is the number of inversions in the permutation, that is the number of pairs $i < j$ such that $\pi(i) > \pi(j)$. Here the order is important since a single value of $j$ can appear multiple times. The order is taken to be the same as the order of the starting positions at the bottom of the walk. An example of walks along the braid $\beta = \sigma_1 \sigma_2^{-1} \sigma_1 \sigma_2^{-1}$ will follow Theorem ~\ref{coloredjones}.

Note that we are reading the braid in two different directions. When writing the braid as a product of generators the braid is read from the top down. When walking along the braid it is read from the bottom up. This is important to get the order of the paths correct in the product of the weights.

Finally a stack of walks is an ordered collection of walks, and the weight of the stack is the product of the weights of walks in the appropriate order. We call this a stack because when considering a picture of these objects it can be thought of as simply stacking the walks on top of each other.

For a stack of walks $W$, we will also consider the local weight at a crossing $j$, denoted $W_j$, which is the product of the weights of the crossing for each walk in the stack. Note that with this notation we can express the weight of $W$ as $(-1)^n(-q)^{\sum_k(|J_k| + \text{inv}(\pi_k))}\prod_{j}W_j$, where the sum ranges over the walks in the stack and $n$ is the number of walks in the stack.

It will be useful to talk about an ordering on the paths that make up a stack. If two paths belong to two different walks at different level of the stack, then the path in the higher walk (that is the walk whose weight is multiplied to the left of the other) is said to be above the other path. If the two paths are in the same walk, then the path which begins to the left of the other path is said to be above the other path.

We will now give a more substantial meaning to the weights of these walks. The letters $a_{j,\pm}$, $b_{j,\pm}$, and $c_{j,\pm}$ will be the operators defined by Huynh and L\^{e} in~\cite{hl}, which will now be repeated here. This will involve introducing a few extra variables. We will eventually make simple substitutions to remove these variables.

Define operators $\hat{x}$ and $\tau_{x}$ and their inverses acting on the ring $\R[x^{\pm1},y^{\pm1},u^{\pm1}]$:
$$\hat{x}f(x,y,\ldots) = xf(x,y,\ldots),\qquad \tau_xf(x,y,\ldots) = f(qx,y,\ldots)$$
Also define $\hat{y}$, $\tau_y$, $\hat{u}$, $\tau_u$, and their inverses similarly.

Now let us define

\begin{tabular}{c c c}
 $a_+ = (\hat{u} - \hat{y} \tau_x^{-1})\tau_y^{-1},$&$b_+ = \hat{u}^2,$&$c_+ = \hat{x}\tau_y^{-2}\tau_u^{-1},$\\
$a_- = (\tau_y -  \hat{x}^{-1})\tau_x^{-1}\tau_u,$&$b_- = \hat{u}^2,$&$c_- = \hat{y}^{-1}\tau_x^{-1}\tau_u$
\end{tabular}

If $P$ is a polynomial in the operators $a_\pm$, $b_\pm$, and $c_\pm$ then we get a polynomial $\E(P) \in \R[z^{\pm1}]$ by having $P$ act on the constant polynomial $1$ and replacing $x$ and $y$ with $z$ and replacing $u$ with $1$. We can get a polynomial $\E_N(P) \in \R$ by making the further substitution in $\E(P)$ replacing $z$ with $q^{(N-1)}$.

We can now define the operators $a_{j,\pm}$, $b_{j,\pm}$, and $c_{j,\pm}$ by replacing the $x$, $y$, $u$ and $z$ in the above definition by $x_j$, $y_j$, $u_j$ and $z_j$. These operators will act on $\Pk=\bigotimes_{j=1}^k \R[x_j,y_j,u_j]$ where $k$ is the number of crossings in $\beta(\gamma)$. It is immediate that any two of these operators with different indecices will commute.

For calculations using these operators, it is useful to observe the relations between operators with the same indecices:

\begin{tabular}{c c c}
$a_+b_+=b_+a_+,$ & $a_+c_+=qc_+a_+,$ & $b_+c_+=q^2c_+b_+$\\
$a_-b_-=q^2b_-a_-,$ & $c_-a_-=qa_-c_-,$ & $c_-b_-=q^2b_-c_-$\\
\end{tabular}

Also we can write a fromula for their evaluation:
\begin{lemma}[HUYNH, L\^{E}]
\label{evaluation}
 $$\E_N(b_+^s c_+^r a_+^d)=q^{r(N-1-d)} \prod_{i=0}^{d-1}(1-q^{N-1-r-i})$$
$$\E_N(b_-^s c_-^r a_-^d)=q^{-r(N-1)}\prod_{i=0}^{d-1}(1-q^{r+i+1-N})$$
\end{lemma}

We can now state the theorem relating the colored Jones polynomial to walks. This is simply a reinterpretation of the main theorem in ~\cite{hl}. The proof will be presented in section ~\ref{qd}.

\begin{theorem}
\label{coloredjones}
Given a braid $\beta(\gamma)$ whose closure is the knot $K$,
\begin{eqnarray*}
J'_K(N) &=& q^{(N-1)(\omega(\beta)-m+1)/2}\sum_{n=0}^\infty \E_N(C^n)\\
&=&q^{(N-1)(\omega(\beta)-m+1)/2} \E_N(\mathscr{S})
\end{eqnarray*}
where the polynomial $C$ is the sum of the weights of walks on $\beta(\gamma)$ with $J \subset \{2, \ldots, m\}$. Furthermore, $\mathscr{S}=\sum_{n=0}^\infty C^n$ is the sum of the weights of the stacks of walks on $\beta(\gamma)$ with $J \subset \{2, \ldots, m\}$.
\end{theorem}

Before presenting an example, there is a simplification we can make to this theorem. It turns out that in general there will be several cancelling terms in this sum. The following lemma states what some of these cancellations are.

\begin{definition}
A simple walk is a walk in which no two paths in the collection traverse the same point on the braid. 
\end{definition}

\begin{lemma}
\label{cancellation}
a) For any nonsimple walk $\beta(\gamma)$, there is another walk whose weight is the negative of the original. The nonsimple walks occur in cancelling pairs.

b) For any stack of walks which traverse the same point on $N$ different levels and has weight $W$, the evaluation $\E_N(W)$ of that weight will be zero.
\end{lemma}

In other words, part a) of this lemma tells us that in Theorem ~\ref{coloredjones} the occurrences of the word ``walks'' may be replaced with ``simple walks''. In later sections, when Theorem ~\ref{coloredjones} is applied, we will assume all walks are simple. Part b) assures us that the sum will be finite and gives us a limit on the stacks of walks we need to consider. The proof of this lemma will follow an example.

\begin{example}
Let $\gamma = ((1,+),(2,-),(1,+),(2,-))$. Thus $\beta = \beta(\gamma) = \sigma_1 \sigma_2^{-1} \sigma_1 \sigma_2^{-1}$ and $K$, the closure of $\beta$, is the figure-eight knot. We also have that $m=3$ and $\omega(\beta) = 0$. The only two walks along $\beta$ which do not start or end at the first strand are presented in the following figure:

\begin{figure}[htbp] %
   \centering
   \includegraphics[width=2in]{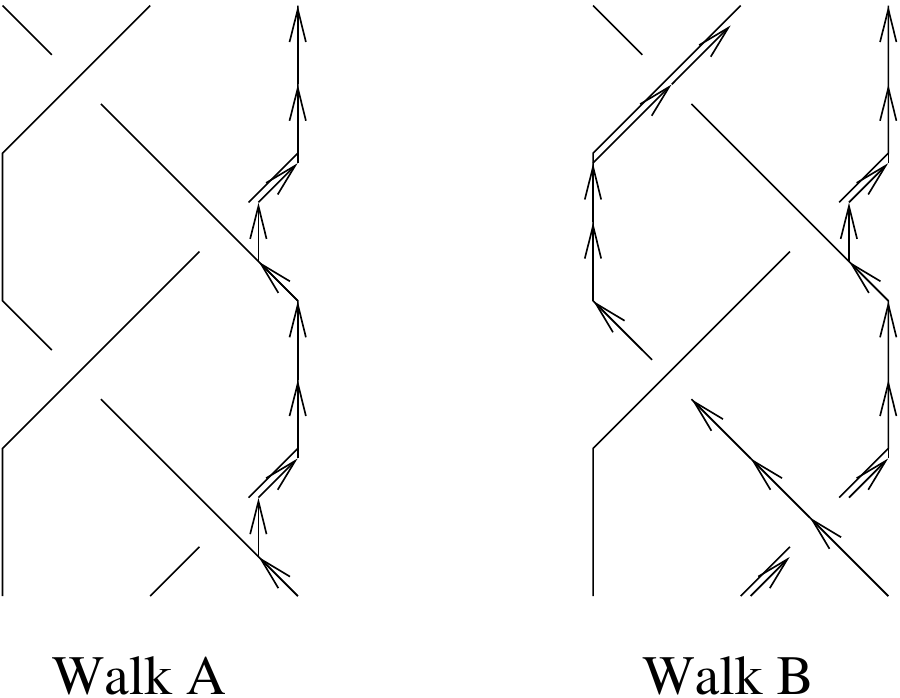} 
   \caption{Two Example Walks}
   \label{walksfig}
\end{figure}

Notice that the walk $A$ consists of a single path from $3$ to $3$ and walk $B$ consists of two paths, one from $2$ to $3$ and one from $3$ to $2$. The only possible walks which do not start at the first strand are paths from $2$ to $2$, of which there are none, paths from $3$ to $3$, the only one being $A$, and walks consisting of two paths, one of which starts at $2$ and the other starts at $3$, and they end at $2$ and $3$ (not necessarily respectively); $B$ is the only walk of this kind. It is an easy exercise to confirm that there are no other such walks.

We will often not distinguish between the weight of a walk and the walk itself. It should be clear from context, for instance $A = q a_{2,-} a_{4,-}$ and $B = q^3 a_{2,-}b_{4,-}c_{1,+}b_{3,+}c_{4,-}$.

Thus by Theorem~\ref{coloredjones} and Lemma~\ref{cancellation}
$$J'_K(N) = q^{(1-N)}\sum_{n=0}^{N-1} \E_N((q a_{2,-}a_{4,-} + q^3 a_{2,-}b_{4,-}c_{1,+}b_{3,+}c_{4,-})^n)$$
The reason the sum stops at $N-1$ is because both $A$ and $B$ traverse the bottom right corner of the braid, and thus any stack with more than $N-1$ levels will evaluate to zero by part b) of Lemma~\ref{cancellation}.

We will use Lemma~\ref{evaluation} to evaluate this sum for given values of $N$.

First for $N=2$,

\begin{eqnarray*}
\E_2(q a_{2,-}a_{4,-}) &=& q(1-q^{-1})^2\\ \\
\E_2(q^3 a_{2,-}b_{4,-}c_{1,+}b_{3,+}c_{4,-}) &=& \E_2(q^3 c_{1,+} a_{2,-} b_{3,+} b_{4,-} c_{4,-})\\
&=& q^3 * q * (1-q^{-1}) * q^{-1}\\ &=& q^3(1-q^{-1})
\end{eqnarray*}

Thus

\begin{eqnarray*}
J'_K(2) &=& q^{-1}(1 + q(1-q^{-1})^2 + q^3(1-q^{-1}))\\
& =& q^2 - q + 1 - q^{-1} + q^{-2}
\end{eqnarray*}

which is the ordinary Jones polynomial for the figure-eight knot.

%For $N=3$,
%$$\E_3(q a_{2,-}a_{4,-}) = q(1-q^{-2})^2$$
%
%$$\E_3(q^3 c_{1,+} a_{2,-} b_{3,+} b_{4,-} c_{4,-}) = q^3 * q^2 * (1-q^{-2}) * q^{-2} = q^3(1-q^{-2})$$
%
%$$\E_3(q^2 a_{2,-}a_{4,-}a_{2,-}a_{4,-}) = q^2(1-q^{-2})^2 (1-q^{-1})^2$$
%
%$$\E_3(q^6 c_{1,+} a_{2,-} b_{3,+} b_{4,-} c_{4,-}c_{1,+} a_{2,-} b_{3,+} b_{4,-} c_{4,-})$$
%$$= \E_3(q^8 c_{1,+}c_{1,+} a_{2,-}a_{2,-} b_{3,+}b_{3,+} b_{4,-}b_{4,-}c_{4,-}c_{4,-})$$
%$$= q^8 * q^4 * (1-q^{-2})(1-q^{-1}) * q^{-4} = q^8(1-q^{-2})(1-q^{-1})$$
%
%$$\E_3(q^4 a_{2,-}a_{4,-} c_{1,+} a_{2,-} b_{3,+} b_{4,-} c_{4,-})$$
%$$= \E_3(q^5 c_{1,+} a_{2,-}a_{2,-} b_{3,+} b_{4,-}c_{4,-}a_{4,-})$$
%$$= q^5 * q^2 * (1-q^{-2})(1-q^{-1}) * q^{-2}(1 - q^{-1}) = q^5(1-q^{-2})(1-q^{-1})^2$$
%
%$$\E_3(q^4 c_{1,+} a_{2,-} b_{3,+} b_{4,-} c_{4,-} a_{2,-}a_{4,-})$$
%$$= \E_3(q^4 c_{1,+} a_{2,-}a_{2,-} b_{3,+} b_{4,-}c_{4,-}a_{4,-})$$
%$$= q^4 * q^2 * (1-q^{-2})(1-q^{-1}) * q^{-2}(1 - q^{-1}) = q^4(1-q^{-2})(1-q^{-1})^2$$

For $N \geq 2$, we need to expand the binomial. Since operators with different subscripts commute, this means that $AB = q BA$ because $a_{4,-}b_{4,-}c_{4,-} = qb_{4,-}c_{4,-}a_{4,-}$. Thus 
$$(A+B)^n = \sum_{k=0}^n \left( \begin{matrix}n\\k\end{matrix}\right)_{q} B^k A^{n-k}$$

$$J'_K(N) = q^{(1-N)}\sum_{n=0}^{N-1} \sum_{k=0}^n \left( \begin{matrix}n\\k\end{matrix}\right)_{q} q^{n+k(k+1)} [\prod_{j=1}^n (1-q^{j-N})] [\prod_{i=1}^{n-k} (1-q^{k+i-N})]$$
where $\left( \begin{matrix}n\\k\end{matrix}\right)_{q}$ is the q-binomial coefficient
$$\left( \begin{matrix}n\\k\end{matrix}\right)_{q} = \prod_{i=0}^{k-1}\frac{1-q^{n-i}}{1-q^{i+1}}$$

\end{example}

\begin{proof}[proof of lemma~\ref{cancellation}]

Part a:  Consider a walk $W$ on the braid $\beta$ where a point on the braid is tranversed more than once, such as the walk in figure ~\ref{badwolf}. There may in general be many such points that are traversed more than once. Consider the highest crossing (equivalently the crossing with the lowest index number) in which two paths in $W$ seperate and call it's index $I$. For example, in figure ~\ref{badwolf}, $I=2$ because the two paths separate at the second crossing from the top, which is crossing $2$. There is another walk, call it $W'$ which passes through the same points as $W$, but at crossing $I$, the two paths that separate take the opposite direction than was taken in $W$. Figure ~\ref{badwolf} shows a pair of walks which differ in this way. Either $W$ or $W'$ has the property that the two paths which pass through crossing $I$ separate so that the one that began to the left of the other one is to the left of the other one immediately after crossing $I$. Denote the walk with this property as $W^{(1)}$ and the other as $W^{(2)}$. Also denote the previously mentioned path that begins to the left of the other one as $X^{(i)}$ and the path that begins to the right of the other one as $Y^{(i)}$, so that $X^{(i)}$ and $Y^{(i)}$ are paths in $W^{(i)}$. In Figure ~\ref{badwolf}, the walk on the left is $W^{(1)}$ and the walk on the right is $W^{(2)}$. Also, the paths consisting of the straight arrows are the $X^{(i)}$'s and the paths consisting of the round arrows are the $Y^{(i)}$'s.

 The claim is that the weight of $W^{(1)}$ is the negative of the weight of $W^{(2)}$, and thus when they are added together in the colored Jones polynomial, they will cancel out. This is clearly a bijection from the set of nonsimple walks to itself, which means we need not consider any walks of this type.

\begin{figure}[htbp] %
   \centering
   \includegraphics[width=1.5in]{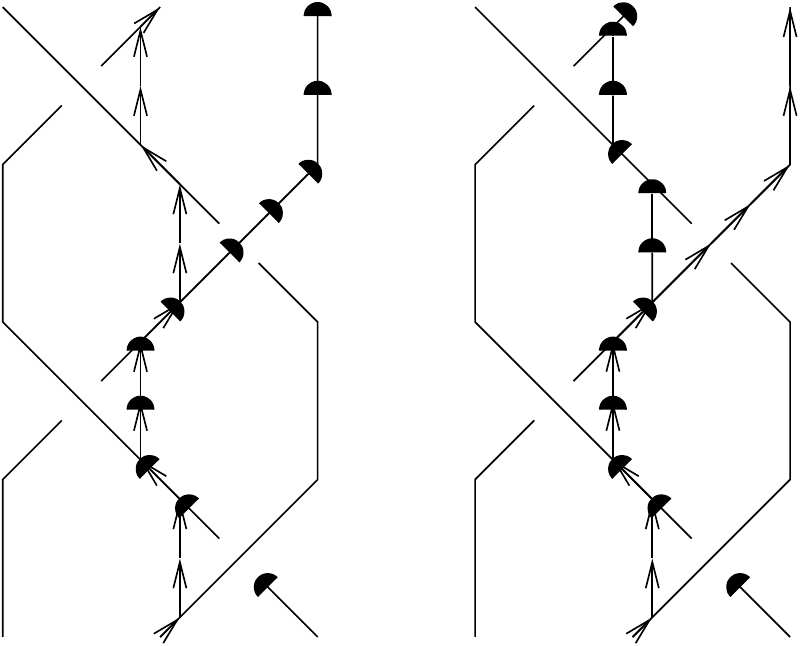} 
	\caption{A Pair of Nonsimple Walk}
	\label{badwolf}
\end{figure}

Denote $W^{(i)}_j$ for the local weight of $W^{(i)}$ at crossing $j$. Observe that at crossing $I$, the walks $W^{(1)}$ and $W^{(2)}$ consists of an $a_{I,\epsilon_I}$ and $c_{I,\epsilon_I}$ in some order. If $\epsilon_I = +$, then $W^{(1)}_I$ is $a_{I,+}c_{I,+}$ and $W^{(2)}_I$ is $c_{I,+}a_{I,+} = q^{-1}a_{I,+}c_{I,+}$. If $\epsilon_I = -$, then $W^{(1)}_I$ is $c_{I,-}a_{I,-}$ and $W^{(2)}_I$ is $a_{I,-}c_{I,-} = q^{-1}c_{I,-}a_{I,-}$. In both cases, $W^{(2)}_I = q^{-1} W^{(1)}_I$.

The paths $X^{(i)}$ and $Y^{(i)}$ may continue to cross above crossing $I$, but since they cannot meet again it must be with $b$'s and $c$'s. If they cross at a positive crossing $j$ with $X^{(i)}$ over $Y^{(i)}$, then $W^{(i)}_j$ will be $c_{j,+}b_{j,+}$. If they cross at a positive crossing $j$ with $Y^{(i)}$ over $X^{(i)}$, then $W^{(i)}_j$ will be $b_{j,+}c_{j,+} = q^2c_{j,+}b_{j,+}$. If they cross at a negative crossing $j$ with $Y^{(i)}$ over $X^{(i)}$, then $W^{(i)}_j$ will be $b_{j,-}c_{j,-}$. And if they cross at a negative crossing $j$ with $X^{(i)}$ over $Y^{(i)}$, then $W^{(i)}_j$ will be $c_{j,-}b_{j,-} = q^2b_{j,-}c_{j,-}$.

If the paths $X^{(i)}$ and $Y^{(i)}$ cross an even number of times above crossing $I$, then the paths $X^{(i)}$ and $Y^{(i)}$ will contribute one inversion in the pemutation associated to $W^{(2)}$ which is not in $W^{(1)}$. %Putting all of this together, the product of the local crossings considered so far gives us that $W_2 = (-q)q^{-1}W_1 = -W_1$.
If the paths $X^{(i)}$ and $Y^{(i)}$ cross an odd number of times above crossing $I$, then the paths $X^{(i)}$ and $Y^{(i)}$ will contribute one inversion in $W^{(1)}$ which is not in $W^{(2)}$. %Again putting this together, the product of the local crossings considered so far gives us $(-q)W_2 = (q^2)q^{-1}W_1$ thus $W_2 = -W_1$.
Assuming no other paths cross $X^{(i)}$ or $Y^{(i)}$ then in both cases we see that $W^{(i)} = -W^{(i)}$, and the proof is complete.

There are, however, several possible cases where another path crosses one of $X_i$ or $Y_i$. We will work through one case, and the other cases can be worked through similarly.

Suppose there is a path $Z$ in $W^{(i)}$ that begins between $X^{(i)}$ and $Y^{(i)}$, and ends between $X^{(i)}$ and $Y^{(i)}$. The only difference in the weights of $W^{(1)}$ and $W^{(2)}$ comes from the crossings between $Z$ and $X^{(i)}$ and $Y^{(i)}$, as well as two additional inversions in the permutation of either $W^{(1)}$ or $W^{(2)}$. The difference between the weights at crossings above $I$ involving $Z$ and one of the special paths will again be $b_\epsilon c_\epsilon$ versus $c_\epsilon  b_\epsilon$ as described earlier. The path $Z$ must cross one of $X^{(1)}$ and $Y^{(1)}$ an even number of times and the other an odd number of times. Thus change in the weights at crossings involving $Z$ will be $q^{\pm 1}$, and the change coming from the inversions involving $Z$ will be $q^{\mp 1}$. Thus in this case we again have $W^{(2)} = -W^{(1)}$.

Part b: Suppose a stack of walks $W$ traverses a point on the braid on $N$ or more different levels. Starting at that point follow along the braid until you reach an over-strand of a crossing. If the top of the braid is reached first, then the same position at the bottom of the braid will have the same starting positions there as there were ending positions, so we can consider the walks to continue from the bottom. Thus we can follow along the braid until we reach an over-strand. The weight at this overstrand will be the product of $N$ or more $c$'s and $a$'s with possibly some additional $b$'s. If there is at least one $a$ and less than $N$ $c$'s, then by lemma ~\ref{evaluation}, the evaluation of the local weight $\E_N(W_j)$ at that crossing will have a factor of $(1-q^0)$ and thus $\E_N(W) = 0$. If the number of $c$'s at this crossing is $N$ or more however, continue along the overstrand until the next overstrand , possibly starting from the bottom if you reach the top of the braid again. If this process continues until you reach the original point where this process started without ever having come across an overstrand with less than $N$ $c$'s taken, then the part of the braid traversed through this process will be a component of the closure of the braid. However, the closure of this braid is a knot, and there is a point on the closure of the braid that could not have been traversed corresponding to the lower left starting position on the braid. Thus, the traversed area could not be a component and thus there must be a crossing with $\E_N(W_j)=0$.

\end{proof}

\section{Quantum Determinants}
\label{qd}

In this section we shall discuss the theorem of Huynh and L\^{e} describing the colored Jones polynomial as the inverse of the quantum determinant of a certain `almost quantum' matrix. All of the details in this section up to the proof of Theorem ~\ref{coloredjones} come from ~\cite{hl}.

A $2 \times 2$ matrix $\left(\begin{matrix}
a & b \\
c & d
\end{matrix}
\right)$ is right quantum if

\begin{eqnarray*}
ac&=&qca\\
bd&=&qdb\\
ad &=& da + qcb -q^{-1}bc
\end{eqnarray*}

An $m\times m$ matrix is right-quantum if all $2\times2$ submatrices of it are right-quantum.

If $A=(a_{ij})$ is right-quantum, then the quantum determinant is 
$$\text{det}_q(A) := \sum_{\pi \in \text{Sym}(m)} (-q)^{\text{inv}(\pi)}a_{\pi 1, 1}a_{\pi 2, 2}\ldots a_{\pi m, m}$$
where inv($\pi$) denotes the number of inversions.

In general, $I-A$, where $I$ is the identity matrix is no longer right-quantum. So define
$$\widetilde{\text{det}}_q(I-A):=1-C$$
where
$$C:=\sum_{\emptyset\neq J \subset \{1,2,\ldots, m\}}(-1)^{|J|-1}\text{det}_q(A_J),$$
where $A_J$ is the $J$ by $J$ submatrix of $A$, which is always right-quantum.

Now that we have defined the almost quantum determinant that we will need, we will define a right quantum matrix from a braid whose closure is a knot. First define matrices which are right quantum:
$$S_+:=\left(\begin{matrix}
a_+ & b_+ \\
c_+ & 0 
       \end{matrix}\right)
\qquad S_-:=\left(\begin{matrix}
0 & c_- \\
b_- & a_-
            \end{matrix}\right)$$

Given a braid $\beta(\gamma)$ defined as in ~\ref{walks}, associate to each $\sigma_{i_j}^{\epsilon_j}$ the matrix which is the identity except for the $2\times2$ minor of rows $i_j$, $i_j + 1$ and columns $i_j$, $i_j + 1$ which is replaced by the matrix $S_{\epsilon_j,j}$.

Here $S_{\pm,j}$ is the same as $S_{\pm}$ with $x$, $y$, $u$ replaced by $x_j$, $y_j$, $u_j$.

The matrix $\rho(\gamma)$ is the product of these matrices.
The matrix $\rho'(\gamma)$ is $\rho(\gamma)$ with the first row and column removed.

\begin{theorem}[Huynh, L\^{e}]
\label{qdetcjp}

$$J'_K(N) = q^{(N-1)(\omega(\beta)-m+1)/2}\E_N \left(\frac{1}{\widetilde{\text{det}}_q(I-q\rho'(\gamma))}\right)$$
 
\end{theorem}

Here $$\frac{1}{\widetilde{\text{det}}_q(I-q\rho'(\gamma))}= \sum_{n=0}^\infty C^n.$$ This sum is finite if the closure of $\beta(\gamma)$ is a knot.

\begin{proof}[proof of Theorem ~\ref{coloredjones}]
 
In order to prove Theorem ~\ref{coloredjones}, we need to show that the polynomial $C$ in Theorem ~\ref{qdetcjp} defined by the quantum determinant, is the same as the polynomial $C$ in Theorem ~\ref{coloredjones} which was defined to be the sum of the weights of the walks along $\beta(\gamma)$.

Step 1: The main idea is that the matrix multiplication corresponds to the choices made during a walk. More explicitly, if $\rho(\gamma) = (M_{i,j})$, then $M_{i,j}$ is the sum of the weights of the paths from $j$ to $i$. We will show this by induction on the length of the braid word.

Base case: $\rho(\emptyset) = I_m$, also $\beta(\emptyset) = e_m$, where $I_m$ is the $m \times m$ identity matrix, and $e_m$ is the identity braid on $m$ strands. It is straight-forward to see that the sum of the weights of paths from $j$ to $i$ along the identity braid is $1$ if $i=j$ and $0$ otherwise.

Inductive step: Now suppose the claim is true for braid words of lenght $k$ and consider a braid word $\beta(\gamma)$ of length $k+1$. If $\gamma_1=(l,+)$, (i.e. the first letter in $\beta(\gamma)$ is $\sigma_l$) then 

$$\rho(\gamma)=(1\oplus \ldots \oplus 1 \oplus
\left(\begin{matrix}
a_{1,+} & b_{1,+} \\
c_{1,+} & 0 
       \end{matrix}\right)
\oplus 1 \oplus \ldots \oplus 1)( \rho(\gamma'))$$

where $\gamma'$ has length $k$.

If $\rho(\gamma) = (M_{i,j})$ and $ \rho(\gamma') = (M'_{i,j})$, then 
$$M_{i,j} = M'_{i,j} \text{   for   } i \neq l, l+1$$
$$M_{l,j} = a_{1,+} M'_{l,j} + b_{1,+} M'_{(l+1),j}$$
$$M_{(l+1),j} = c_{1,+} M'_{l,j}$$

Let us now compare this with the paths along $\beta(\gamma)$.

\begin{figure}[htbp] %
   \centering
   \includegraphics[width=1.5in]{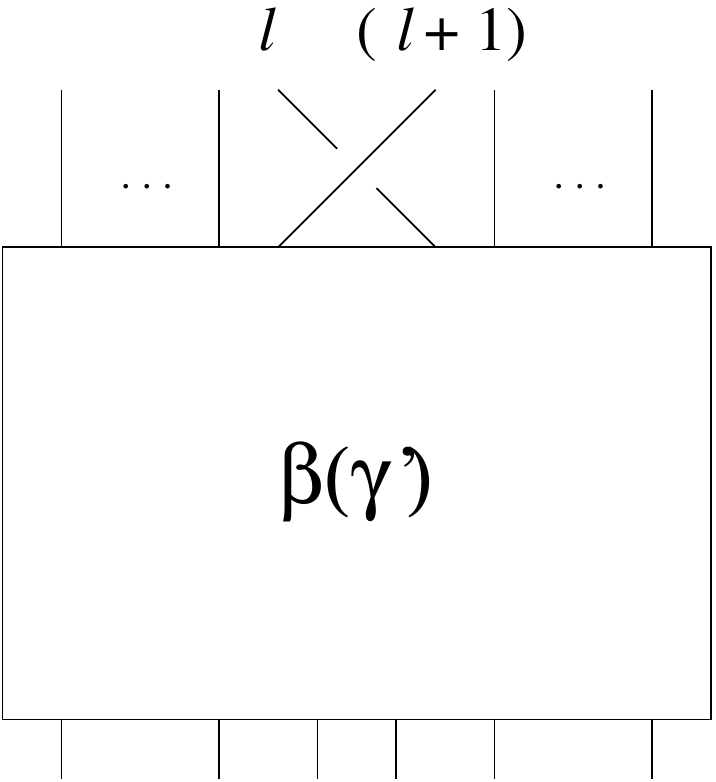} 
   \caption{}
\end{figure}

By induction, $M'_{i,j}$ is the sum of the weights of the paths along $\beta(\gamma')$ from $j$ to $i$, so the sum of the weights along $\beta(\gamma)$ from $j$ to $i$ when $i\neq l,(l+1)$ is also $M'_{i,j}$. The paths from $j$ to $l$ come in two types: those that walk along $\beta(\gamma')$ from $j$ to $l$ and then jump down at crossing $1$, and those that walk along $\beta(\gamma')$ from $j$ to $(l+1)$ and follow along the lower strand of crossing $1$. Thus the sum of the weights of these paths is $a_{1,+} M'_{l,j} + b_{1,+} M'_{(l+1),j}$. Finally, the paths from $j$ to $(l+1)$ consists of paths along $\beta(\gamma')$ from $j$ to $l$ and then following along the upper strand of crossing $1$. Thus the sum of the weights of these paths is $c_{1,+} M'_{l,j}$. This completes step 1.

Step 2: The rest of the proof is simply following through the definitions of the weights of walks and the inverse of the quantum determinant.

The polynomial $C$ from Theorem \ref{qdetcjp} is the sum
$$C = \sum_{\emptyset\neq J \subset \{1,2,\ldots, m\}}(-1)^{|J|-1}\text{det}_q(\rho'(\gamma)_J).$$
Since $\rho'(\gamma)$ is just $\rho(\gamma)$ with the first row and column removed, we can write $C$ as
$$C = \sum_{\emptyset\neq J \subset \{2,\ldots, m\}}(-1)^{|J|-1}\text{det}_q(\rho(\gamma)_J).$$
For a particular subset $J$, the expression $(-1)^{|J|-1}\text{det}_q(\rho(\gamma)_J)$ is precisely the sum of the weights of the walks for the given $J$. Thus $C$ is the sum of the weights of all walks for all $J$.

\end{proof}

\section{Positive Braids: Proof of Theorem 1}
\label{posbraids}

Suppose that $\beta(\gamma)$ is a positive braid, meaning that $\epsilon_j = +$ for all $j$. To prove Theorem ~\ref{positive} we will show that, in the sum $\sum_n \E_N(C^n)$, every monomial with $n>0$ has degree at least $N$.

Let us consider lowest terms in the evaluation of a stack of walks $W$. In order to apply lemma ~\ref{evaluation} we need to rearrange the order of the product. Since the operators corresponding to different crossings commute, we can rearrange the product and evaluate the weight at each crossing.

Define $W_j$ to be the product of the local weights at crossing $j$. Let $A_j$, $B_j$, and $C_j$ be the number of $a$'s, $b$'s, and $c$'s respectively in the local weight $W_j$ at crossing $j$. Also let ${A_B}_j$ be the number of pairs of an $a$ and a $b$ in the local weight at crossing $j$ such that the $a$ is to the left of the $b$, in other words, the number of ``commutations'' that would need to be made to arrange the letters so that all $b$'s appear to the left of all $a$'s. Define, similarly, notation for all combinations of $A$, $B$, and $C$.

For a stack of walks $W$, by Lemma ~\ref{evaluation} the lowest degree will be:
$$\sum_k (|J_k|+\text{inv}(\pi_k)) + \sum_j w_j$$
where $w_j = {A_C}_j - 2{C_B}_j -C_j A_j +(N-1)C_j$

%Using a similar reasoning as the agument at the beginning of Lemma ~\ref{3}, we will add or subtract $1$ at each crossing where two paths cross as in the following table. (The table is reduced because there are only positive crossings.)

The goal now is to find a useful lower bound for this sum. Firstly, now that we have an explicit sum for the minimum degree, we will modify the terms in the sum without changing the total value. At each crossing, we will add or subtract $1$ to the term $w_j$ every time two paths cross each other. If the path originally on the left is above the other path, then add $1$; if the path originally on the right is above the other path subtract $1$. The first column of the following table shows how paths might cross each other and the result of adding or subtracting $1$ from each crossing of paths that occurs at that crossing in the braid. The second column shows the situation where two paths may come together (this can only happen if the two paths are on different level of the stack) and then seperate, which may or may not count as a crossing of paths. If the paths seperate without crossing then we may add $1$ at one of the crossings and subtract $1$ at the other, so that resultant sum is not changed. If the paths do cross from this situation, then we will add or subtract $1$ at only one of the crossings.

\begin{center}

\begin{tabular}{|c|c|c|}
\hline & \includegraphics[width=.3in]{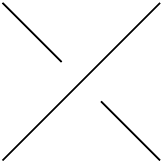}& \includegraphics[width=.2in]{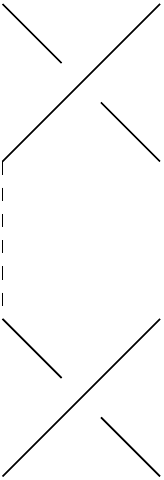} \\ \hline
L over R (top)  &  $+{C_B}_j$ & $+{C_A}_j$\\ \hline
R over L (bottom)& $-{B_C}_j$ & $-{B_A}_j$\\ \hline
\end{tabular}

\end{center}

By adding or subtracting $1$ in this way throughout, we can describe how the total sum will change. The original sum was  $\sum_{j} w_{j}$. Call the new sum $\sum_j v_j$, where $v_j$ is the new term coming from crossing $j$ after all of the additions and subtractions are applied. If you restrict to only paths on the same walk, say the $k$-th walk, then the result of all of the additions and subtractions will be a change by $\text{inv}(\pi_k)$. If you restrict to two different walks and consider only the $\pm 1$'s coming which occurs between paths on these two different levels, then all of these additions and subtractions will result in no change to the total sum. This can be seen as follows:

Each of the walks themselves can be thought of as braids. The closure of these braids are links. If we stack the two links coming from these two walks the way the walks are stacked, then we can see that the additions and subtractions are just the calculation of the sum of the linking numbers of different components of these links, where the different components come from the different links. This is obviously $0$ since one link is entirely above the other.

Thus 
$$\sum_{j} v_j = \sum_k \text{inv}(\pi_k) + \sum_{j} w_j$$

To every crossing, we have added:
$C_B -B_C -B_A +C_A$

Thus
\begin{eqnarray*}
v_j  &=& (N-1)C_j - C_j B_j - {B_A}_j \\
\end{eqnarray*}

Now, because $\sum_j B_j = \sum_j C_j$, we can define

\begin{eqnarray*}
u_j  &:=& (N-1)B_j - C_j B_j - {B_A}_j \\
 &\geq & (N-1)B_j - C_j B_j - B_j A_j \\
     &=& (N-1-C_j-A_j)B_j \\
 &\geq & 0
\end{eqnarray*}

And we get $$\sum_{j} u_j = \sum_{j} v_j$$

There will necessarily be a crossing which has some number of $b$'s and no $a$'s or $c$'s. Call this crossing $\iota$ and we then have
$$u_\iota \geq (N-1-C_\iota-A_\iota)B_\iota \geq N-1$$

Finally we get that the lowest degree in $\E_N(W)$ is
\begin{eqnarray*}
&& \sum_k (|J_k|+\text{inv}(\pi_k)) + \sum_j w_j\\
&=&\sum_k |J_k| + \sum_j v_j\\
&=&\sum_k |J_k| + \sum_j u_j\\
&\geq & 1 + u_\iota \\
&\geq & N
\end{eqnarray*}


\begin{thebibliography}{99}

\bibitem{ck} A. Champanerkar, and I. Kofman, {\it On the Tail of Jones Polynomials of closed braids with a full twist},

%\bibitem{dl} O. Dasbach, and X.-S. Lin, {\it On the head and the tail of the colored Jones polynomial}, Compositio Math., Vol 142 (2006), No. 5, pp. 1332-1342.

\bibitem{hl} V. Huynh, and TTQ. L\^{e}, {\it On the colored Jones polynomial and the Kashaev invariant}, Fundam. Prikl. Mat. \textbf{11} (2005), no. 5,57-78.

\bibitem{j} V. Jones, {\it Hecke algebra representations of braid group and link polynomials}, Ann. Math. \textbf{126} (1987), 335-388.

\bibitem{ltw} X.-S. Lin, F. Tian, and Z. Wang, {\it Burau representation and random walk on a string link}, Pacif. J. Math., \textbf{182} (1998), 289-302.

\bibitem{lw} X.-S. Lin, Z. Wang, {\it Random walk on knot diagrams, colored Jones polynomial and Ihara-Selberg zeta function}, Proceedings of the Birmanfest, AMS/IP Studies in Advanced Mathematics, vol. 24, 2001, pp. 107--121.



\end{thebibliography}
\end{document}